\documentclass[review]{elsarticle}

\usepackage{lineno,hyperref}
\modulolinenumbers[5]

\def\R{\mathbb{R}}
\usepackage{latexsym,amsfonts}
\usepackage[active]{srcltx} 
\newtheorem{theorem}{Theorem}[section]

\newtheorem{notation}[theorem]{Notation}

\newtheorem{corollary}[theorem]{Corollary}
\newtheorem{lemma}[theorem]{Lemma}

\newtheorem{remark}[theorem]{Remark}

\begin{document}

\begin{frontmatter}

\title{Infinite solutions having a prescribed number of nodes for a p-Laplacian problem}




\author{Jing Zeng\corref{mycorrespondingauthor}}
\cortext[mycorrespondingauthor]{Corresponding author}
\ead{zengjing@fjnu.edu.cn}
\address{School of Mathematics and Computer Sciences, Fujian
Normal University, Fuzhou, P.R.China}


\begin{abstract}
In this paper, we are concern with the multiplicity of solutions for a p-Laplacian problem. A weaker super-quadratic assumptions is required on the nonlinearity.
Under the weaker condition we give a new proof for the infinite solutions having a prescribed number of nodes to the problem. It turns out that the weaker condition on nonlinearity suffices to guarantee the infinitely many solutions.
 At the same time, a global characterization of the critical values of the nodal radial solutions are given.

\end{abstract}


\end{frontmatter}

\linenumbers

\section{Introduction}
Many paper are concerned with the existence and multiplicity of radial solutions and non-radial solutions of the semilinear equation
$$-\Delta u+u=f(x,u),\quad \quad u\in H^{1}(\mathbb{R}^N). \eqno{(1.1)}$$
The equation originates from various problems in the field of physics and mathematical physics. (1.1) is called a Euclidean field equation
in cosmology \cite{CGM1978}. And nonlinear Klein-Gordon or Schr\"odinger equations when one is looking for certain types of solitary waves \cite{BL1983,BL19832,S1977}. More general, (1.1) can be explained as the case of $p=2$ in the more general problem
$$-\Delta_p u+|u|^{p-2}u=f(x,u),\quad \quad u\in W^{1, p}(\mathbb{R}^N). \eqno{(1.2)}$$

Since (1.1) is invariant under rotations it is natural to search for spherically symmetric solutions. For the radial solutions of (1.1) is proved by Bartsch-Willem \cite{BW1993}, Liu-Wang \cite{LW2004} and Li \cite{L1990}. The existence of non-radial solutions of (1.1) or (1.2) seems to have been open for a long time \cite{DN1986}. The
non-radial solutions of (1.1) were proved by Bartsch-Willem \cite{BW19932} and Liu-Wang \cite{LW2004}.

For the p-Laplacian equation $$-\Delta_p u=f(x,u),\ u\in W^{1, p}_0(\Omega),\eqno{(1.3)}$$ where $\Omega$ is a bounded domain in $\R^N$, Dinca-Jebelean-Mawhin \cite{DJM2001} obtained the existence results under Dirichlet boundary condition. Bartscha-Liu \cite{BL2004} proved the existence of four solutions for  equation (1.3), that is, a pair of subsolution and supersolution, a positive and a negative solution, in addition a sign changing solution. Bonanno-Candito \cite{BC2003} established the existence of three solutions to the Neumann boundary condition of problem (1.2).

\smallskip

We require the following assumptions on the nonlinearity $f(x, u)$:
\smallskip

{\bf $(f_1)$} $f(x, 0)=0, f(x, t)=o(|t|^{p-2}t),$ as $|t|\rightarrow0,$ uniformly in $x$.

\smallskip

{\bf $(f_2)$}  $f\in C(\R^N, \R)$ and there exist $C>0$, $q\in (p, p^*)$ such that
$$|f(x, t)|\leq C(1+|t|^{q-1}),$$
where $p^*=Np/(N-p)$ if $N>p$, and $p^*=\infty$ if $N\leq p$.

\smallskip

{\bf $(f_3)$} $\displaystyle\lim_{|t|\rightarrow\infty}\frac{F(x, t)}{|t|^p}=+\infty$, where $F(x, t)=\displaystyle\int^t_0f(x, s)ds$.

\smallskip

{\bf $(f_4)$} There exists $R>0$ such that, for any $x$, $\displaystyle\frac{f(x, t)}{|t|^{p-2}t}$ is increasing in $t\geq R$,
and decreasing in $t\leq -R$.

\smallskip

\begin{remark} \label{th-1}
The condition $(f_3)$ is a consequence of the following condition: $$\lim_{|t|\rightarrow\infty}\frac{f(x, t)}{|t|^{p-2}t}=+\infty.$$
In the case $p=2$, $(f_3)$ characterizes the problem (1.2) as superlinear at infinity. It is a extension of  a much natural super-quadratic conditions {\bf $((SQ)$} condition$)$,

\smallskip

{\bf $(SQ)$} $\displaystyle\lim_{|t|\rightarrow\infty}\frac{F(x, t)}{t^2}=\infty.$

\end{remark}

\smallskip

\begin{remark} \label{th-2}
$(SQ)$ condition is weaker than the famous Ambrosetti-Rabinowitz
growth condition $((AR)$ condition$)$. Since the seminal work of Ambrosetti-Rabinowitz \cite{AR1973}, the (AR) condition is most frequent appeared
in the superlinear elliptic boundary value problem.

{\bf $(AR)$} There exist $\mu>p$ and $R>0$ such that
$$ 0<\mu F(x, t)\leq f(x, t) t, \ \ \hbox{for}\ x\in D \ \hbox{and}\ |t|\geq R.$$
It is important not only in establishing the mountain-pass geometry of the functional but also in obtaining the bounds of (PS) sequences.
 In fact, {\bf $(AR)$} condition implies that for some $C>0$,
$$F(x, t)\geq C|t|^{\mu}, \ \mu>p .$$

\end{remark}

In recent years there were some articles, such as \cite{WZ2003, MS2008}, trying to drop the $(AR)$ condition in the studying of the
superlinear problems. For equation (1.1), Liu-Wang \cite{LW2004} first posed the {\bf $(SQ)$} condition
to get the bounds of minimizing sequence on Nehari manifold. Furthmore under coercive condition of potential function $V(x)$, they proved the
existence of three solutions of equation $-\Delta u+V(x)u=f(x,u)\ (u\in H^1(\R^N))$, one positive, one negative and one sign-changing solution.
Li-Wang-Zeng \cite{LWZ2006} gave a natural generalization of the results in \cite{LW2004} to two noncompact cases. Miyagaki-Souto \cite{MS2008} established the existence of nontrivial solution of (1.1) by combining some arguments used by
Struwe-Tarantello \cite{ST1998}. Liu \cite{L2010} obtained the existence and multiplicity results for equations (1.3), and considered the Cerami sequences of the Euler-Lagrange functional.

\smallskip

The main result of this paper is
\smallskip

\begin{theorem} \label{th-3} Under assumptions
$(f_1)-(f_4)$, for every integer $k>0$, there exist a pair $u^{+}_k$ and $u_{k}^-$ of radial solutions of (1.2) with $u_{k}^-(0)<0<u_{k}^+(0)$, having exactly
$k$ nodes $0<\rho_1^{\pm}<\cdots<\rho_{k}^{\pm}<\infty$.
\end{theorem}

\smallskip

Here a node $\rho>0$ is such that $u(\rho)=0$.

\smallskip

\begin{theorem} \label{th-8}
Under assumptions $(f_1)-(f_4)$, $f(x, u)$ is odd in $u$, there exist infinitely many non-radial nodal solutions of (1.2).
\end{theorem}

\smallskip

\begin{remark} \label{th-9}
It is also possible to replace the oddness of $f(x, u)$ by other conditions, we refer the reader to the work of Jones-K\"upper \cite{JK1986}.
\end{remark}

\smallskip

If we further assume that:

{\bf $(f_5)$}  $f\in C(\R^N, \R)$ and for some $C>0$,
$$|f'_u(x, t)|\leq C(1+|t|^{q-2}),$$
where $q=p^*$ if $N\geq3$ and $q\in (p, p^*)$ if $N=2.$

\smallskip

\begin{corollary} \label{th-10}
Assume $N=4$ or $N\geq6$, assume $(f_2)-(f_5)$ hold, $f$ is odd in $u$, then equation (1.2) has an unbounded sequence of non-radial
sign-changing solutions.
\end{corollary}

\smallskip
In the present paper, we give a new proof for the infinite solutions having a prescribed number of nodes to the problem (1.2), and the results are got under the weaker (SQ) conditions. It turns out that the (SQ) condition on $f(x, u)$ suffices to guarantee infinitely many solutions.
Our method is spirted by the work of \cite{LW2004}, and our theorems generalize the results in \cite{LW2004} to the case of $p\neq2$. At the same time, a global characterization of the critical values of the nodal radial solutions are given.

\section{Preliminaries}
In this section, we give some notations and some preliminaries lemmas, which will be adopted in the proof of Theorems.

Solutions of (1.2) is correspond to the critical points of the functional
$$J(u):=\int_{\R^N}\frac{1}{p}|\nabla u|^p+\frac{1}{p}|u|^p-F(x,u), \ u\in W^{1,p}(\R^N),$$ where
$W^{1,p}(\R^N)$ is endowed with the norm
$\displaystyle\|u\|=\bigg(\int_{\R^N}(|\nabla u|^p+|u|^p)\bigg)^\frac{1}{p}$.

\smallskip

\begin{notation} \label{th-4} We define the Nehari manifold
$$\mathcal{N}_1=\{u\in X_1: u\neq 0, \langle J'(u), u\rangle=0\},$$ where $X_1:=\{u\in W^{1,p}(\R^N):u(x)=u(|x|)\}$.
And $$\mathcal{N}_2=\{u\in X_2: u\neq 0, \langle J'(u), u\rangle=0\},$$ where $X_2:=W^{1,p}(\R^N)$.
For $0\leq\rho<\sigma\leq\infty$, define
$$\Omega(\rho, \sigma):= int\{x\in\R^N: \rho\leq|x|\leq\sigma\},$$  $$X_{\rho, \sigma}:=\{u\in W^{1,p}(\Omega(\rho, \sigma)):u(x)=u(|x|)\},$$
$$\mathcal{N}_{\rho, \sigma}=\{u\in X_{\rho, \sigma}: u\neq 0, \langle J'(u), u\rangle=0\}.$$
Define $u(x)=0$ for $x \notin \Omega(\rho, \sigma)$ if $u\in X_1$.Obviously $X_{\rho, \sigma}\subset X_1$ and $\mathcal{N}_{\rho, \sigma}\subset \mathcal{N}_1$.

Fix $k$, define
$$\mathcal{N}^{+}_k=\{u\in X_1: \hbox{there\ exist}\ 0=\rho_0<\rho_1<\cdots<\rho_{k}<\rho_{k+1}=\infty \  \hbox{such\ that}$$
$$(-1)^ju|_{\Omega(\rho_j, \rho_{j+1})}\geq0 \  \hbox{and}\ u|_{\Omega(\rho_j, \rho_{j+1})}\in \mathcal{N}_{_{\rho_j, \rho_{j+1}}} \hbox{for}\ j=0, \cdots,k \}.$$

Defined on $[0, \infty)\times\R$,
$$f^+(r, u)=\left\{
\begin{array}{ll}
f(r, u),\qquad\  \hbox{if} \ u\geq0,\\
-f(r, -u),\quad \hbox{if} \ u<0, \\
\end{array} \right.$$
and $F^+(r, u)=\displaystyle\int_0^uf^+(r, s)ds$,
$$J^+(u):=\int_{\R^N}\frac{1}{p}|\nabla
u|^p+\frac{1}{p}|u|^p-F^+(x,u).$$
Similarly we can define
$$f^-(r, u)=\left\{
\begin{array}{ll}
f(r, u),\qquad\ \hbox{if} \ u\leq0,\\
-f(r, -u),\quad \hbox{if} \ u>0,\\
\end{array} \right.$$
and $F^-(r, u)$, $J^-(u)$.
\end{notation}

The letters $C$ will always denote various universal constants.

\smallskip

\begin{lemma} \label{th-5} Under assumptions $(f_1)-(f_4)$, equation
$$-\Delta_p u+|u|^{p-2}u=f(x,u),\quad \quad u\in X_{\rho, \sigma}, \eqno{(2.1)}$$
has a weak solution $u$ such that
$$J(u)=\max_{t>0}J(tu)=\inf_{v\in X_{\rho, \sigma}\backslash\{0\}}\max_{t>0}J(tv)>0.$$
\end{lemma}
\noindent{\bf Proof.} By the assumptions $(f_1)$ and $(f_2)$, $J$ has a strict local minimum at 0. For any $u\neq0$, $ J(tu)\rightarrow-\infty$ as $t\rightarrow\infty$. Thus $$c:=\inf_{v\in X_{\rho, \sigma}\backslash\{0\}}\max_{t>0} J(tv)>J(0)=0\eqno{(2.2)}$$ is well-defined.

Let $(u_n)$ be a minimizing sequence of $c$ such that
$$ J(u_n)=\max_{t>0} J(tu_n)\rightarrow c$$ as $n\rightarrow\infty$.

First we want to prove that $(u_n)$ is bounded. If not, consider $v_n:=u_n/\|u_n\|,$ then $\|v_n\|=1$.
By passing to a subsequence, we may assume $v_n\rightarrow v$ weakly in $X_{\rho, \sigma}$ and strongly in $L^r(X_{\rho, \sigma})$ for any $r\in [p, p^*].$
Note that $(f_1)$ and $(f_2)$ implies $\int_{X_{\rho, \sigma}}F(x, u)$ is weakly continuous on $X_{\rho, \sigma}$.

If $v\neq0$, we have $$\frac{c+o(1)}{\|u_n\|^p}=\frac{1}{p}-\int_{X_{\rho, \sigma}}\frac{F(x, u_n)}{u_n^p}v_n^p. $$
By (2.2), $$\frac{1}{p}>\int_{X_{\rho, \sigma}}\frac{F(x, u_n)}{u_n^p}v_n^p.$$
Then by $(f_3)$ and Fadou's lemma, passing a limit on the both sides, then
$$\frac{1}{p}>\int_{X_{\rho, \sigma}}\frac{F(x, u_n(x))}{u_n^p}v^p=\infty.$$ It gives a contradiction.

 If $v=0$, fixing an $R>\sqrt[p]{pc}$, by $\|v_n\|=1$, we have
$$ J(u_n)\geq J(Rv_n)=\frac{1}{p}R^p-\int_{X_{\rho, \sigma}}F(x, Rv_n).$$
$J(u_n)$ converges towards $c$, but $R^p/p-\int_{X_{\rho, \sigma}}F(x, Rv_n)$ tends to $R^p/p>c$, a contradiction. Thus $(u_n)$ is bounded.

Assume $u_n$ weakly converges to $u$. As $n\rightarrow\infty$, then
$$\int_{X_{\rho, \sigma}}u_nf(x, u_n)\rightarrow\int_{X_{\rho, \sigma}}uf(x, u).$$
Since, for some $\alpha>0$, $\|u_n\|^p>\alpha$, and
$$\|u_n\|^p=\int_{X_{\rho, \sigma}}u_nf(x, u_n),$$ so $u\neq0$.

There is $s>0$ such that $ J(su)=\max_{t>0} J(tu)$. Then
$$ J(su)\leq\liminf_{n\rightarrow\infty} J(su_n)\leq\liminf_{n\rightarrow\infty} J(u_n)=c.$$

$(f_4)$ implies that $\max_{t>0} J(tu)$ is achieved at only one point $t=s$.
It is also the unique one such that $\langle J'(tu), u\rangle=0$.

Next we claim that $su$ is a critical point of $J$. Without loss of generality, we assume $s=1$. If $u$ is not a critical point, there is $v\in C_0^{\infty}(\Omega)$ such that $\langle J'(u), v\rangle=-2$. There is $\varepsilon_0>0$ such that for $|t-1|+|\varepsilon|\leq\varepsilon_0$, $\langle J'(tu+\varepsilon v), v\rangle\leq-1.$

For $\varepsilon>0$ small, let $t_{\varepsilon}>0$ be the unique number such that $$\max J(tu+sv)= J(t_{\varepsilon}u+\varepsilon v).$$ Then $t_{\varepsilon}\rightarrow1$ as $\varepsilon\rightarrow0$.

For $\varepsilon$ small such that $|t_{\varepsilon}-1|+\varepsilon\leq\varepsilon_0$, then $ J(t_{\varepsilon}u+\varepsilon v)\geq c$, but by the assumption that $\langle J'(tu+\varepsilon v), v\rangle\leq-1,$
$$ J(t_{\varepsilon}u+\varepsilon v)= J(t_{\varepsilon}u)+\int^1_0\langle J'(t_{\varepsilon}u+s\varepsilon v), \varepsilon v\rangle ds\leq c-\varepsilon<c.$$ It is a contradiction.

\smallskip

\begin{lemma} \label{th-6}
Under assumptions $(f_2)-(f_5)$, $f$ is odd in $u$, then equation (2.1) has infinitely many pairs of solutions.
\end{lemma}
\noindent{\bf Proof.} It is clear that the solutions occur in pairs due to the oddness of $f(x, u)$. Under the assumptions, any critical point of $J$ restricted on $\mathcal{N}_2$ is a critical point of $J$ in $X_2$. To verify the (PS) condition it suffices to show any (PS) sequence is bounded. This is similar to the proof of Lemma \ref{th-5}. We omit the details. If (PS) condition is satisfied on $\mathcal{N}_2$ then the standard Ljusternik-Schnirelmann theory gives rise to an unbounded sequence of critical values of $J$, see the details in \cite{Ra1986}.

\medskip

\section{Proof of Theorems}
In this section, we prove Theorem \ref{th-3} and Theorem \ref{th-8}.

\noindent {\bf Proof of Theorem \ref{th-3}:}\ First by Lemma \ref{th-5}, the infimum
$$c^+(\rho, \sigma):=\inf_{\mathcal{N}_{\rho, \sigma}}J^+.$$
is achieved. Since $|u|$ is also a minimizer, we may assume the minimizer $u$ is a positive solution of the problem
$$-\Delta_p u+|u|^{p-2}u=f(x,u),\quad \quad u\in X_{\rho, \sigma}. \eqno{(3.1)}$$
Similarly, the infimum $$c^-(\rho, \sigma):=\inf_{\mathcal{N}_{\rho, \sigma}}J^-,$$
is also achieved by negative minimizers which are negative solutions of (3.1).

Then we work on the Nehari manifold $\mathcal{N}_{k}^+$, and construct a $u^{+}_k\in \mathcal{N}^{+}_k$ such that  $$\displaystyle c^{+}_k:=\inf_{\mathcal{N}^{+}_k}J$$ is achieved by some $u^{+}_k$, which gives the desired solutions in Theorem \ref{th-3}.

Let $(u_n)$ be a minimizing sequence of $c^{+}_k$. As the same arguments in the proof of Lemma \ref{th-5}, $(u_n)$ is bounded.

Since $u_n\in \mathcal{N}^{+}_k$, there exist $0=\rho_0^n<\rho_1^n<\cdots<\rho_{k}^n<\rho^n_{k+1}=\infty$ such that
$(-1)^ju_n|_{\Omega(\rho_j^n, \rho_{j+1}^n)}\geq0$ and $u_n|_{\Omega(\rho_j^n, \rho_{j+1}^n)}\in \mathcal{N}_{_{\rho_j^n, \rho_{j+1}^n}}$ for
$j=0, \cdots,k$.

Note that
$$\|u_n|_{\Omega(\rho_j^n, \rho_{j+1}^n)}\|^p=\int_{\Omega(\rho_j^n, \rho_{j+1}^n)}u_nf(r, u_n).$$
By $(f_1)-(f_2)$, 0 is a strict local minimizer of $ J$, thus there is a $\delta>0$ such that $\|u\|\geq\delta$ for $u\in\mathcal{N}_{\rho_j^n, \rho_{j+1}^n}$. Fix $q\in (p, p^*)$, for any $\varepsilon>0$, there is a constant $C>0$ such that
$$\int_{\Omega(\rho_j^n, \rho_{j+1}^n)}u_nf(r, u_n)\leq
\varepsilon\int_{\Omega(\rho_j^n, \rho_{j+1}^n)}|u_n|^p+C\int_{\Omega(\rho_j^n, \rho_{j+1}^n)}|u_n|^q,$$
where $q\in(p, p^*)$.
Therefore, by choosing $\varepsilon>0$ small we can find a $C>0$ such that
$$\delta^p\leq\|u_n|_{\Omega(\rho_j^n, \rho_{j+1}^n)}\|^p\leq C\int_{\Omega(\rho_j^n, \rho_{j+1}^n)}|u_n|^q.\eqno{(3.2)}$$
Using (3.2), in a similar way as in \cite{BW1993}, one sees that $(\rho_{k+1}^n)_n$ is bounded away from $\infty$, $(\rho_{j+1}^n-\rho_{j}^n)_n$
is bounded away from 0 for each $j$,
and there are $0=\rho_0<\rho_1<\cdots<\rho_{k}<\rho_{k+1}=\infty$ such that $\rho_j^n\rightarrow\rho_j$ as $n\rightarrow\infty$, for $j=1, \cdots, k$.

Along a subsequence of $(n)$, we may assume that $u_n\rightarrow u$ weakly in $X_1$, strongly in $L^r(X_1)$ for any
$r\in [p, p^*]$. It follows that $u_n|_{\Omega(\rho_j^n, \rho_{j+1}^n)}\rightarrow u|_{\Omega(\rho_j, \rho_{j+1})}$ weakly in
$X_1$, strongly in $L^r(X_1)$
($r\in [p, p^*)$).
And $(-1)^ju|_{\Omega(\rho_j, \rho_{j+1})}\geq0$, for $u\in\mathcal{N}_{\rho_j^n, \rho_{j+1}^n}$.

Letting $n\rightarrow\infty$ in (3.2),
it implies that $u|_{\Omega(\rho_j, \rho_{j+1})}\neq0$. Thus we can choose an $\alpha_j>0$ such that
$\alpha_ju|_{\Omega(\rho_j, \rho_{j+1})}\in\mathcal{N}_{(\rho_j, \rho_{j+1})}$ for $j=1, \cdots, k.$ Define
$$u_k^+:=\sum^k_{j=0}\alpha_ju|_{\Omega(\rho_j, \rho_{j+1})}.$$
By the definition of $u_k^+$, it can observe that $u_k^+\in\mathcal{N}_k^+$.

Next we want to show

1) $c_k^+$ is archived by $u_k^+$, that is, $ J(u_k^+)=c_k^+$,

2) $u_k^+$ is a radial function having nodes
$0<\rho_1<\cdots<\rho_{k}<\infty$,

3) $u_k^+$ is a solution of (1.2).

The weak convergence of $u_n|_{\Omega_{(\rho_j^n, \rho^n_{j+1})}}$ in $X_1$ and strong convergence in $L^r(X_1)(p<r<p^*)$ imply
$$c_k^+\leq J(u_k^+)=\sum_{j=0}^k J(\alpha_ju|_{\Omega(\rho_j,\rho_{j+1})})\leq
\sum_{j=0}^k\liminf_{n\rightarrow\infty} J(\alpha_ju_n|_{\Omega(\rho_j^n, \rho^n_{j+1})}). \eqno{(3.3)}$$
And $$\sum_{j=0}^k\liminf_{n\rightarrow\infty} J(u_n|_{\Omega(\rho_j^n, \rho^n_{j+1})})=
\liminf_{n\rightarrow\infty} J(u_n)=c_k^+.\eqno{(3.4)}$$
So $ J(u_k^+)=c_k^+$.

Then the equality in (3.3) implies that $\alpha_ju|_{\Omega(\rho_j^n, \rho^n_{j+1})}$ is a minimizer of $$\displaystyle\inf_{\mathcal{N}_{\rho_j^n, \rho^n_{j+1}}\cap P^{+}} J^{+},\ \hbox{if}\ j\ \hbox{is}\ \hbox{even},$$ and a minimizer of $$\displaystyle\inf_{\mathcal{N}_{\rho_j^n, \rho^n_{j+1}}\cap P^{-}} J^{-},\ \hbox{if}\ j\ \hbox{is}\ \hbox{odd},$$
where $P^{\pm}:=\{u\in X_1:\pm u\geq0\}$. At the same time,
$\alpha_ju|_{\Omega(\rho_j, \rho_{j+1})}$ is a minimizer of $\inf_{\mathcal{N}_{\rho_j, \rho_{j+1}}} J^{\pm}.$
For $j$ even, $\alpha_ju|_{\Omega(\rho_j^n, \rho^n_{j+1})}$ is a positive solution of (1.2), and for $j$ odd, $\alpha_ju|_{\Omega(\rho_j^n, \rho^n_{j+1})}$ is a negative solution.
Then the strong maximum principle implies that $u_k^+(0)>0$, $(-1)^ju_k^+(x)>0$, for $\rho_j<|x|<\rho_{j+1}$ ($j=0, 1, \cdots, k$), and $$(-1)^j\lim_{|x|\uparrow\rho_j}\frac{\partial u_k^+(x)}{\partial |x|}>0,\ \
(-1)^j\lim_{|x|\downarrow\rho_j}\frac{\partial u_k^+(x)}{\partial |x|}>0, \ \ \hbox{for}\ j=1, \cdots, k.$$  So $u_k^+$ has exactly $k$ nodes.

In order to prove $u_k^+$ is a solution of (1.2), for simplicity we assume $\alpha_j=1$ for all $j$. If $u_k^+$ is
not a critical point of $J$, then there is a $ \varphi\in C_0^{\infty}(\R^N)$ such that $$\langle J'(u_k^+),  \varphi\rangle=-2.$$ Observe that there is an $\tau>0$ such that
if $|s_j-1|\leq\tau$ ($ j=0, \cdots, k$) and $0\leq\varepsilon\leq\tau$ then the function $$g(s, \varepsilon):=\displaystyle\sum_{j=0}^ks_ju|_{\Omega(\rho_j, \rho_{j+1})}+\varepsilon \varphi,$$ where $s=(s_1, \cdots, s_k)$, has exactly $k$ nodes
$0<\rho_1(s, \varepsilon)<\cdots<\rho_k(s, \varepsilon)<\infty.$ And
$\rho_j(s, \varepsilon)$ is continuous in $(s, \varepsilon)\in D\times[0, \tau]$, where $D:=\{(s_1, \cdots, s_k)\in \R^k:|s_j-1|\leq\tau\}$, and
$$\bigg\langle J'(g(s, \varepsilon)),  \varphi\bigg\rangle<-1.\eqno{(3.5)}$$

In order to deduce a contradiction, we set for $s\in D$,
$$g_1(s)=\sum_{i=0}^ks_iu|_{\Omega(\rho_i, \rho_{i+1})}+\tau\eta(s) \varphi,$$
where $\eta(s): D\rightarrow[0, 1]$ ($s=(s_1, \cdots, s_k)$) is a cut-off function such that
$$\eta(s_1, \cdots, s_k)=\left\{
\begin{array}{ll}
1,\quad \hbox{if} \ |s_i-1|\leq\tau/4 \ \ \hbox{for all}\ i,\\
0,\quad \hbox{if} \ |s_i-1|\geq\tau/2 \ \ \hbox{for at least one}\ i.\\
\end{array} \right.$$
Then for each $s\in D$, $g_1(s)\in C(D, X)$, and $g_1(s)$ has exactly $k$ nodes $0<\rho_1(s)<\cdots<\rho_k(s)<\infty$,
where $\rho_j(s)$ is continuous.

Further, we define for $j=1, \cdots, k$,
$$h_j(s):=\bigg\langle J'(g_1(s))\bigg|_{\Omega(\rho_j(s), \rho_{j+1}(s))}, g_1(s)\bigg|_{\Omega(\rho_j(s), \rho_{j+1}(s))}\bigg\rangle.$$
And define $h: D\rightarrow\R^k$ as $h(s):=(h_1(s), \cdots, h_k(s))$. Then $h(s)\in C(D, \R^k).$

For a fixed $j$, if $|s_j-1|=\tau$ then $\eta(s)=0$ and $\rho_i(s)=\rho_i$ for all $i=1, \cdots, k$. So by the definition of $g_1(s)$,

$$h_j(s)=\langle J'(s_ju)|_{\Omega(\rho_j, \rho_{j+1})}, s_ju|_{\Omega(\rho_j, \rho_{j+1})}\rangle
=\left\{
\begin{array}{ll}
>0,\quad \hbox{if} \ s_j=1-\tau,\\
<0,\quad \hbox{if} \ s_j=1+\tau.\\
\end{array} \right.$$
Therefore, the degree $deg(h, int(D), 0)$ is well defined and $deg(h, int(D), 0)=(-1)^k$. Thus there is an $s\in int(D)$ such that $h(s)=0$, that is, $g_1(s)\in \mathcal{N}_k^+$.

It is obviouly
$$ J(g_1(s))\geq c_k^+.\eqno{(3.6)}$$

On the other hand, by (3.5),
$$ J(g_1(s))= J(\sum_{j=0}^ks_ju|_{\Omega(\rho_j, \rho_{j+1})})+\int^1_0\bigg\langle J'(\sum_{j=0}^ks_ju|_{\Omega(\rho_j, \rho_{j+1})}+\theta\tau\eta(s) \varphi), \tau\eta(s) \varphi\bigg\rangle d\theta$$
$$\leq J(\sum_{j=0}^ks_ju|_{\Omega(\rho_j, \rho_{j+1})})-\tau\eta(s).\qquad\qquad\qquad\qquad\quad\qquad\qquad\ \ \ \ $$

If $|s_j-1|\leq\tau/2$ for each $j$, then by (3.4)
$$ J(g_1(s))< J(\sum_{j=0}^ks_ju|_{\Omega(\rho_j, \rho_{j+1})})\leq\sum_{j=0}^k J(u|_{\Omega(\rho_j, \rho_{j+1})})=c_k^+,\eqno{(3.7)}$$
which contradicts (3.6).

If $|s_j-1|>\tau/2$ for at least one $j$, by (3.4)
$$ J(g_1(s))\leq J(\sum_{j=0}^ks_ju|_{\Omega(\rho_j, \rho_{j+1})})<\sum_{j=0}^k J(u|_{\Omega(\rho_j, \rho_{j+1})})=c_k^+,\eqno{(3.8)}$$
A contradiction with (3.6) too. The proof is finished.

\smallskip

\noindent {\bf Proof of Theorem \ref{th-8}:} Using a result of Lions \cite{PL1982}, it is possible to fine a subspace $E$ of $X_2$ consisting of functions which are not radial and such that the inclusion $E\hookrightarrow L^s$ is compact for $p<s<p^*$, see the detail in Theorem IV.1 of \cite{PL1982} or the proof of Theorem 2.1 in \cite{BW19932}. Then follow the same steps in Lemma \ref{th-5}, and combine Lemma \ref{th-6} to get the infinitely many non-radial nodal solutions of (1.2).

\section*{Acknowledgements}
The author is supported by  China Postdoctoral Science Foundation (2014M551830) and the Nonlinear Analysis Innovation Team
(IRTL1206) funded by Fujian  Normal  University.

\end{document}